\documentclass[12pt,oneside]{article}
\usepackage{amsmath,amssymb,amsfonts,amsthm}
\usepackage{amscd}
\usepackage{graphicx}
\usepackage[all]{xy}
\usepackage{multirow}
\usepackage{rotating}
\usepackage{hyperref}
\usepackage{hyperref}
\hypersetup{
    colorlinks=true,
    linkcolor=blue,
    filecolor=magenta,
    urlcolor=cyan,
    citecolor=blue}
\usepackage{lscape}
\newcommand{\T}{{\cal T}}

\newcommand{\To}{\longrightarrow}

\newcommand{\tm}{\T M}

\newcommand {\cp}{\mathfrak{X}(\pi (M))}

\newcommand{\Ima}{\text{Im}}

\def\Section#1{\vspace{30truept}\addtocounter{section}{1}\setcounter{thm}{0}
\setcounter{equation}{0}{\noindent\Large\bf
    \arabic{section}.~~#1}\par \vspace{12pt}}

\newtheorem{thm}{Theorem}[section]
\newtheorem{cor}[thm]{Corollary}
\newtheorem{lem}[thm]{Lemma}
\newtheorem{prop}[thm]{Proposition}
\newtheorem{defn}[thm]{Definition}

\newtheorem{rem}[thm]{Remark}

\numberwithin{equation}{section}

\begin{document}
\title{ Tripathi Connection in Finsler Geometry}
\author{\bf{A. Soleiman$^{1,2}$ and Ebtsam H. Taha$^{3}$}}
\date{}
\maketitle                     
\vspace{-1.16cm}
\begin{center}
{$^{1\,}$ Department of  Mathematics, Collage of Sciences and Arts - Gurayat, Jouf University, Kingdom of Saudi Arabia}\\
{ $^{2\,}$ Department of Mathematics, Faculty of Sciences, Benha University, Benha,  Egypt}\\
{ $^{3}$ Department of Mathematics, Faculty of Science, Cairo University, Giza, Egypt.}
\end{center}
\vspace{-0.8cm}

\begin{center}
E-mails: amr.hassan@fsc.bu.edu.eg, amrsoleiman@yahoo.com\\ {\hspace{1.3cm}}
ebtsam.taha@sci.cu.edu.eg, ebtsam.h.taha@hotmail.com
\end{center}
\begin{center}
\textit{\textbf{Dedicated to Professor Nabil L. Youssef on the occasion of his 74th Birthday}}
\end{center}
 \maketitle
\smallskip
\noindent{\bf Abstract.} Adopting the pullback formalism,  a new  linear connection in Finsler geometry has been introduced  and investigated.
Such connection unifies all formerly known Finsler connections and some other connections not introduced so far. Also, our connection is a Finslerian version  of the Tripathi connection introduced in Riemannian geometry. The existence and uniqueness of such connection is proved intrinsically. An explicit intrinsic  expression relating this connection to Cartan connection is obtained.  Some generalized Finsler connections are constructed from Tripathi Finsler connection, by applying the  ${P}^{1}$-process and ${C}$-process introduced by Matsumoto.  Finally, under certain conditions, many special  Finsler connections are given.




\vspace{30truept}\centerline{\Large\bf{Introduction}}\vspace{12pt}
\par
In Riemannian geometry, the Levi-Civita connection is the unique torsion-free connection that preserves the Riemannian metric. Hayden introduced a metric connection with torsion  \cite{Hayden}. Folland in \cite{Folland}, with the help of a $1$-form, explored  a symmetric connection that is non-metric. Yano in \cite{Yano.2} investigated  a certain type of Hayden connection which is known as  semi-symmetric metric connection. Then, semi-symmetric non-metric connection had been studied in cf. \cite{{Agashe}, Sengupta}. A further extension of semi-symmetric connections is the notion of quarter-symmetric connection \cite{Yano.1} which includes  Ricci quarter-symmetric connection. Recently, a generalized quarter-symmetric connection has been introduced in \cite{gen quart}.  Fortunately, in \cite{New}, Tripathi defined a new connection  which included all these connections and more as particular cases.\\

\vspace*{-0.3 cm}
Finsler geometry is a natural generalization of Riemannain geometry \cite{r42}. The connection theory in the context of Finsler geometry has been enormously developed in both the Klein-Grifone approach and the pullback approach, see for example \cite{r22,  r74,{Ma.book2}, r93, new1, r86, r94}. There are four celebrated  Finsler connections, namely, Cartan, Berwald, Chern-Rund and Hashiguchi connections. Some types of the above-mentioned Riemannian connections have been extended to the Finsler context \cite{ hv, r42,holonomy, h-rec.}. In \cite{amr20}, Tripathi connection has been introduced locally in the Finsler framework. \\

\vspace*{-0.3 cm}
Based on the above discussion, this paper is devoted to a further development of the theory of connections in Finsler geometry. Our geometric treatment avoids the use of coordinate indices. We provide a Finslerain extension of Tripathi connection (Theorem \ref{th.v1})  which not only includes all the aforementioned connections 
 but much more. Successfully, we derive the relation between some geometric objects associated with such connection and  the corresponding ones of Cartan connection.  Followed by an investigation of its spray,  nonlinear connection, torsion and curvature tensors along with Bianchi identities. Then,  in \S 3, using the ${P}^{1}$-process and ${C}$-process defined in \cite{Ma.book2}, we give a generalized version of the four celebrated Finsler connections. Finally,  in \S 4 as a consequence 26 particular regular Finsler connections are given.

\Section{Connections, Sprays and Finsler metrics}
Here we recall the necessary material for better understanding  the present paper.\\
 For an $n$ dimensional smooth manifold $M$, consider the  tangent bundle  $\pi: T M\longrightarrow M$ and its differential $d\pi: TT M\longrightarrow TM$. The vertical bundle $V(TM)$ of $TM$ is just $ \ker (d\pi)$. Let us denote the pullback bundle of the tangent bundle by $\pi^{-1}(T M)$.
 Let $\mathfrak{F}(TM)$ denote the algebra of $C^\infty$ functions on $TM$ and $\cp$ the $\mathfrak{F}(TM)$-module of differentiable sections of   $\pi^{-1}(T M)$. The elements of $\mathfrak{X}(\pi (M))$ will be called $\pi$-vector fields and denoted by barred letters $\overline{X}$. \\
 \vspace*{-0.4cm}
\par
Now, we recall the short exact sequence of vector  bundle morphisms \cite{r21, {new1}}
$$0\longrightarrow
 \pi^{-1}(TM)\stackrel{\gamma}\longrightarrow T(\T M)\stackrel{\rho}\longrightarrow
\pi^{-1}(TM)\longrightarrow 0 ,\vspace{-0.1cm}$$ where $\T M$ is the slit tangent bundle, $\gamma$ is the natural injection and $\rho := (\pi_{TM}, \pi)$.

The  tangent structure of $TM$ or the vertical endomorphism is the endomorphism $J: T\T M \mapsto T \T M$ defined by $J=\gamma \circ \rho$. Note that  $J^{2}=0, \,\,[J,J ]=0$ and $\ker J = \Ima \, J = V(\T M) $. The Liouville vector field is the vector field given by $ \, \mathcal{C}:=\gamma\, \overline{\eta}, $ where  $\overline{\eta}(u)=(u,u)$ for all $u\in \T M$.
\begin{flushleft}
\textbf{Linear  connections on the pullback bundle $\pi^{-1}(TM)$} \cite{new1, {r94}}.
\end{flushleft}
\vspace{-0.2cm}
Let $D$ be  a linear connection on the pullback bundle $\pi^{-1}(TM)$. 
 The map \vspace{-0.15cm}  \[K:T \T M\longrightarrow
\pi^{-1}(TM):X\longmapsto D_X \overline{\eta} .\]  is called the
connection map of $D$.  The connection $D$ is regular if at each $u\in \T M$, we have the  splitting 
$$ T_u (\T M)=V_u (\T M)\oplus H_u (\T M),$$
where $H_u (\T M):= \{ X \in T_u
(\T M) \,|\,K(X)=0 \}$ is the horizontal space at $u$ .

When $M$ is equipped with a regular connection $D$,  the
   maps $
 \gamma,\,\, \rho |_{H(\T M)}$ and $K |_{V(\T M)}$
 are vector bundle isomorphisms. In this case, $\beta:=(\rho |_{H(\T M)})^{-1}$
  is called the horizontal map of $D$.
\begin{defn}
 The  torsion tensor $\mathbb{T}$ of a regular connection $D$ on
$\pi^{-1}(TM)$ with horizontal map $\beta$ has the following two counterparts:
\begin{description}
\item[(a)] $(h)$h-torsion tensor $Q (\overline{X},\overline{Y})=\mathbb{T}(\beta \overline{X},\beta \overline{Y})$,
\item[(b)] $(h)$hv-torsion tensor $ T(\overline{X},\overline{Y})=\mathbb{T}(\gamma
\overline{X},\beta \overline{Y})$.
\end{description}
   \vspace{-0.2cm}
The curvature tensor $\mathbb{K}$ of  $D$  has the following three counterparts:
\begin{description}
\item[(c)] $h$-curvature tensor $R(\overline{X},\overline{Y})\overline{Z}=\mathbb{K}(\beta
\overline{X},\beta \overline{Y})\overline{Z}$,
\item[(d)] $hv$-curvature tensor $ {P}(\overline{X},\overline{Y})\overline{Z}=\mathbb{K}(\beta
\overline{X},\gamma \overline{Y})\overline{Z}$,
\item[(e)] $v$-curvature tensor ${S}(\overline{X},\overline{Y})\overline{Z}=\mathbb{K}(\gamma
\overline{X},\gamma \overline{Y})\overline{Z}$.
\end{description}
\end{defn}
Consequently, the contracted curvature tensors of a connection $D$ (denoted by $\widehat{{R}}$, $\widehat{ {P}}$ and $\widehat{ {S}}$)  are given,  respectively, by
$$\widehat{ {R}}(\overline{X},\overline{Y})={ {R}}(\overline{X},\overline{Y})\overline{\eta},\quad
\widehat{ {P}}(\overline{X},\overline{Y})={
{P}}(\overline{X},\overline{Y})\overline{\eta},\quad \widehat{
{S}}(\overline{X},\overline{Y})={
{S}}(\overline{X},\overline{Y})\overline{\eta}$$
and are called the $v(h)$-torsion, $v(hv)$-torsion and $v(v)$-torsion, respectively.
\begin{flushleft}
\textbf{Geometry of sprays and Finsler metrics} \cite{r21,r22, {new1}, r92, r94}.
\end{flushleft}
\vspace{-0.2cm}
 A vector field $X$ on $TM$ is said to be spray on $M$  if  $\rho\circ X = \overline{\eta}$ and $[\mathcal{C},X]= X$. Each spray induces  canonically a nonlinear connection $\Gamma := [J,S]$,  which is  homogeneous (i.e., $[\mathcal{C}, \Gamma]=0$). The existence of $\Gamma$ is equivalent to the existence of an $n$-dimensional distribution $H : u \in \tm \longrightarrow H_u\in T_u(\tm)$  supplementary to the vertical distribution; it is called the horizontal distribution. The corresponding horizontal and vertical projectors are given, respectively,  by \vspace{-0.2cm}
\begin{equation}
  \label{projectors}
    h:=\frac{1}{2}  (Id_{TM} + \Gamma), \,\,\,\,\,\,            v:=\frac{1}{2}(Id_{TM} - \Gamma).
\end{equation}
\begin{defn}
A smooth Finsler structure on $M$ is a map  $L: TM \To [0,\infty) $~such~that{\em:}
 \begin{description}
    \item[(a)] $L $ is  $C^{\infty}$ on  $\T M$,  $C^{1}$ on $TM$,
    \item[(b)]$L$ is positively homogeneous of degree $1$ in the directional argument $y$, that is
    $\mathcal{L}_{\mathcal{C}} L=L$, where $\mathcal{L}$ is the Lie derivative,
\item[(c)] The Hilbert $2$-form
    $\Omega:=\frac{1}{2}\,dd_{J}L^{2}$  has a maximal rank.
 \end{description}
The Finsler metric $g$ induced by $\,L\,$ on $\pi^{-1}(TM)$  is defined as follows $$g(\rho X,\rho Y):=\Omega(JX,Y), \ \forall  X, Y \in
    \mathfrak{X}(TM).$$
 \end{defn}
Unlike Riemannian geometry which
has one canonical linear connection on $M$, Finsler geometry admits at least four linear connections on $\pi^{-1}(TM)$:
 Cartan, Chern-Rund, Hashiguchi  and  Berwald connections \cite{r94}. It should be noted that these four connections
are  regular with  $T(\overline{X},
\overline{\eta})=0$.

\par
Every Finsler structure determines uniquely a
spray $G$, called the geodesic or canonical spray \cite{r27}. For the geodesic spray $G$ there exists a unique homogenous nonlinear  connection $\Gamma =
[J,G]$, called the Barthel or canonical connection  associated with the Finsler structure $L$. 
\begin{defn}\label{sem.1} Let $(M,L)$ be a Finsler manifold and  $D$ be a regular connection on
$\pi^{-1}(TM)$ with horizontal map $\beta$. Then, the vector field defined by  $S=\beta \, \overline{\eta}$ is called the spray associated with $D$.
In addition, the vector valued $1$-form  $\Gamma_{D}:=2\,\beta\circ\rho-I$ is a nonlinear connection, called associated with $D$.
\end{defn}
\begin{lem}\label{eqv.}  Let $(M,L)$ be a Finsler manifold and  ${D}$ be a regular connection
on $\pi^{-1}(TM)$ whose connection map is $K$ and whose horizontal
map is $\beta$. A necessary and sufficient condition for the $(1,1)$-type tensor defined by ${\Gamma} =\beta\circ\rho - \gamma\circ K$ to be a nonlinear
    connection on $M$ is that the $(h)$hv-torsion  ${T}$ of ${D}$ satisfies
    ${T}( \overline{X},\overline{\eta})=0$. Thereby, ${\Gamma}$
    coincides with the nonlinear connection associated with $D$. That is,\\ ${\Gamma}=\Gamma_{D}=
    2\beta\circ\rho-I$. Consequently, $h_{\Gamma}=h_{D}=\beta\circ\rho$
    and  $\,v_{\Gamma}=v_{D}=\gamma\circ K$.
\end{lem}

\vspace{-0.7cm}

\Section{Tripathi Finsler connection}
Let us start with the following definition before going to the main result of the paper (Theorem \ref{th.v1}).

\begin{defn}\label{quarter}
A regular connection $\overline{D}$ on $\pi^{-1}(TM)$ is said to be quarter-symmetric if there exist a scalar 1-form $u$ and a vector
  1-form $\varphi$ on $\pi^{-1}(TM)$ such that  the $(h)$h-torsion $\overline{Q}$ of $\overline{D}$ satisfies\,{\em:}
\begin{equation}\label{Q}
\overline{Q}(\overline{X},\overline{Y})=
  u(\overline{Y}) \varphi(\overline{X})-u(\overline{X}) \varphi(\overline{Y}) \quad \forall \,  \overline{X}, \, \overline{Y}\in  \cp.
  \end{equation}
 The $1$-forms $u$ and $\varphi$ are called the  quarter-symmetric forms of $\overline{D}$.

 \smallskip
 In particular, if $\varphi=id_{\pi^{-1}(TM)}$, then $\overline{D}$ is called semi-symmetric. Moreover, if $\varphi=0$ or $u=0$, then $\overline{D}$ is called symmetric. Further, if the pullback bundle is equipped with a Finsler structure and $\varphi=\emph{Ric}_{o}$, where $\emph{Ric}_{o}$ is the  horizontal  Ricci $(1,1)$-type tensor of the Cartan connection, then  $\overline{D}$ is said to be Ricci quarter-symmetric.
 \end{defn}

\begin{thm} \label{th.v1}  Let $(M,L)$ be a Finsler manifold. For  given  functions $f_{1}, f_{2} \in \mathfrak{F}(TM)$,   scalar $1$-forms $A, B, u$ and a vector  $1$-form $\varphi$ on $\pi^{-1}(TM)$, there
exists a unique regular connection $\overline{D}(f_{1}, f_{2}, A, B,u,\varphi)$, or simply $\overline{D}$, on $\pi^{-1}(TM)$ such
that
\begin{description}
  \item[(I)] The horizontal covariant derivative of $g$ with respect to $\overline{D}$ has the form\,{\em:}
  $$(\overline{D}_{{\bar{\beta}}\, \overline{X}} \,g)(\overline{Y},\overline{Z})=2f_{1}A(\overline{X})\,g(\overline{Y},\overline{Z})
  +f_{2}\{B(\overline{Y})\, g(\overline{Z},\overline{X})+B(\overline{Z}) \,g(\overline{X},\overline{Y})\},$$
  \item[(II)] The metric $g$ is  $\overline{D}$-vertically parallel, that is  $\overline{D}_{\gamma \overline{X}}\, \,g=0$,

  \item[(III)] $\overline{D}$ is quarter-symmetric with quarter-symmetric forms $u$ and $\varphi$,
  \item[(IV)] The $(h)$hv-torsion $\overline{T}$ of $\overline{D}$ satisfies   $g(\overline{T}(\overline{X},\overline{Y}),\overline{Z}) =g(\overline{T}(\overline{X},\overline{Z}),\overline{Y})$.
\end{description}
This connection will be named Tripathi Finsler connection and denoted by $GC\Gamma$.
\end{thm}

\begin{proof}Suppose that $(M,L)$ admits some regular connection $\overline{D}$ satisfying (I) - (IV).  We prove the the uniqueness of $\overline{D}$.

 Making use of  (II),  (IV) and Lemma \ref{eqv.},  the associated  nonlinear
connection $\Gamma_{\overline{D}}$ is given by \vspace{-0.2cm}
$$\Gamma_{\overline{D}}=\bar{\beta}\circ\rho -\gamma\circ \bar{K} = \bar{h} - \bar{v} .$$ Here,  $\bar{h}, \,\bar{v}, \, \bar{\beta},\, \bar{K}$ are the horizontal projector, vertical projector, horizontal map  and connection map associated with $\overline{D}$. From (II), (IV) and applying the Christoffel trick, we obtain for all  $\overline{X}, \overline{Y}, \overline{Z}\in \mathfrak{X}(\pi (M))$
\begin{eqnarray}\label{eq.t1}
2g(\overline{D} _{\gamma \overline{X}} \overline{Y}, \overline{Z})&=&\gamma \overline{X}\cdot g(\overline{Y},\overline{Z})
 +g(\overline{Y},\rho[\bar{\beta} \,\overline{Z}, \gamma \overline{X}])+g(\overline{Z},\rho[\gamma \overline{X},\bar{\beta}\, \overline{Y}]).
\end{eqnarray}
Since the difference between two nonlinear connections is a semi-basic vector $1$-form on $TM$ \cite{r21, r94}, we get \vspace{-0.3cm}
\begin{equation}\label{diff}
 \bar{\beta} \,\overline{X}={\beta} \overline{X}+\gamma \overline{X}_{t}, \,\,\text{for some } \overline{X}_{t} \in \cp,
\end{equation}
where $\beta$ is the horizontal map of the Cartan
connection $\nabla$.  As the vertical endomorphism $J$  satisfies $\rho \circ J=0$, we have
\begin{equation}\label{eq.t2}
  \rho[\gamma \overline{X},\bar{\beta}\, \overline{Y}]=\rho[\gamma \overline{X},{\beta} \overline{Y}].
\end{equation}
Considering \cite[Theorem 4(a)]{r92} together with  (\ref{eq.t1}) and (\ref{eq.t2}), it follows that the vertical counterpart of $\overline{D} $ and $\nabla$ coincides, that is \vspace{-0.2cm}
\begin{equation}\label{11}
 \overline{D} _{\gamma \overline{X}} \overline{Y}=\nabla  _{\gamma \overline{X}} \overline{Y}.\vspace{-0.1cm}
\end{equation}
\par
Now, by (I) and (III), we conclude that
\begin{eqnarray}\label{eq.t3}
2\,g(\overline{D} _{\bar{\beta} \,\overline{X}} \overline{Y}, \overline{Z})&=&\bar{\beta}\, \overline{X}\cdot g(\overline{Y},\overline{Z})
+\bar{\beta}\, \overline{Y}\cdot g(\overline{Z},\overline{X})-\bar{\beta}\, \overline{Z}\cdot g(\overline{X},\overline{Y})-2\,f_{2}\,B(\overline{Z})\, g(\overline{X},\overline{Y})\nonumber\\
&&-u(\overline{Z}) \{ g(\overline{X}, \varphi(\overline{Y}))+ g(\overline{Y}, \varphi(\overline{X}))\} +g(\overline{Z},\rho[\bar{\beta}\, \overline{X},\bar{\beta}\, \overline{Y}])\nonumber\\
&&+u(\overline{Y})\,\{g(\overline{X}, \varphi(\overline{Z})) + g(\overline{Z}, \varphi(\overline{X})) \}+g(\overline{Y},\rho[\bar{\beta}\, \overline{Z},\bar{\beta}\, \overline{X}]) \nonumber\\
&&
+u(\overline{X}) \{g(\overline{Y},\varphi(\overline{Z}))- g(\overline{Z},\varphi(\overline{Y}))\}-g(\overline{X},\rho[\bar{\beta}\, \overline{Y},\bar{\beta}\, \overline{Z}]) \nonumber\\
&& -2\,f_{1}\{A(\overline{X})\,g(\overline{Y},\overline{Z}) +A(\overline{Y})\,g(\overline{Z},\overline{X}) -A(\overline{Z})\,g(\overline{X},\overline{Y}) \}.
\end{eqnarray}
Formula (\ref{diff}) gives rise to
\begin{equation*}
  \rho[\bar{\beta}\, \overline{X},\bar{\beta}\, \overline{Y}]=\rho[\beta \overline{X},{\beta} \overline{Y}]+
  \rho[\beta \overline{X},{\gamma} \overline{Y}_{t}]+\rho[\gamma \overline{X}_{t},{\beta} \overline{Y}].
\end{equation*}
Then, by  \cite[Theorem 4(b) and Theorem 6]{r92},  (\ref{eq.t3}) becomes \vspace{-0.1cm}
\begin{eqnarray}\label{eq.t4}
g(\overline{D} _{\bar{\beta}\, \overline{X}} \overline{Y}, \overline{Z})&=&g(\nabla _{{\beta} \overline{X}} \overline{Y}, \overline{Z})
+g(D^{\circ} _{{\gamma} \overline{X}_{t}} \overline{Y}, \overline{Z})- f_{2}\,B(\overline{Z})\,g(\overline{X},\overline{Y})\nonumber \\
&&+ \textbf{T}( \overline{X}_{t},\overline{Y}, \overline{Z})  + \textbf{T}( \overline{Y}_{t},\overline{Z}, \overline{X})-\textbf{T}( \overline{Z}_{t},\overline{X}, \overline{Y}) \nonumber\\
&&
-u(\overline{Z})\,g(\varphi_{1}(\overline{X}), \overline{Y})+u(\overline{Y})\,g(\varphi_{1}(\overline{X}), \overline{Z})-u(\overline{X})\,g(\varphi_{2}(\overline{Y}), \overline{Z}) \nonumber\\
&&- f_{1}\{A(\overline{X})\,g(\overline{Y},\overline{Z})+A(\overline{Y})\,g(\overline{Z},\overline{X})-A(\overline{Z})\,g(\overline{X},\overline{Y})\},
\end{eqnarray}
where  $\textbf{T}$ is the Cartan  tensor defined by $\textbf{T}( \overline{X},\overline{Y},\overline{Z}):=g(T( \overline{X},\overline{Y}), \overline{Z})$, $g(\varphi_{1}(\overline{X}), \overline{Y})$ and $g(\varphi_{2}(\overline{X}), \overline{Y})$ are the symmetric and antisymmetric parts of
$g(\varphi(\overline{X}), \overline{Y})$,  respectively.\\
Setting  $\overline{X}=\overline{Y}=\overline{\eta}$ in (\ref{eq.t4}),  using the property that $\textbf{T}$ is indicatory  and  $\bar{K} \circ \bar{\beta}=K \circ \beta=0$,  we get
\begin{eqnarray}\label{B}
\overline{\eta}_{t}&=&f_{1}\{2\,A(\overline{\eta})\,\overline{\eta}-L^{2}\,\overline{a}\}+f_{2}\,L^{2}\,\overline{b} +L\, \ell(\varphi_{1}(\overline{\eta}))\,\overline{u}-u(\overline{\eta})\,(\varphi_{1}-\varphi_{2})(\overline{\eta}),
\end{eqnarray}
where $\ell := L^{-1}\, i_{\overline{\eta}}g,\,\,\,g(\overline{a},\overline{X}):=A(\overline{X})$,  $\,\,g(\overline{b},\overline{X}):=B(\overline{X}),\,\,\,g(\overline{u},\overline{X}):=u(\overline{X})$.\\
 Set $\overline{Y}=\overline{\eta}$ again in  (\ref{eq.t4}) and consider (\ref{B}), we obtain
 \begin{eqnarray}\label{B1}
   \overline{X}_{t}&=&f_{1}\{A(\overline{X})\,\overline{\eta}
   +A(\overline{\eta})\,\overline{X}-L\,\ell(\overline{X})\, \overline{a} +L^{2}\, T(\overline{a},\overline{X})\}\nonumber\\
   &&     + f_{2}\{L\,\ell(\overline{X})\, \overline{b} -L^{2}\, T(\overline{b},\overline{X})\} -u(\overline{\eta})\,\{\varphi_{1}(\overline{X}) + T((\varphi_{2}-\varphi_{1})(\overline{\eta}), \overline{X})\} \nonumber\\
   &&
   +L\,\{\ell(\varphi_{1}(\overline{X}))\,\overline{u}-\ell(\varphi_{1}(\overline{\eta}))\, T(\overline{u},\overline{X})\}  +  u(\overline{X})\varphi_{2}(\overline{\eta}) .
 \end{eqnarray}
 Therefore, the Cartan tensor satisfies
   \begin{eqnarray}
   \textbf{T}(\overline{X}_{t},\overline{Y}, \overline{Z})&=&f_{1}\,\{A(\overline{\eta})\,\textbf{T}(\overline{X},\overline{Y}, \overline{Z})
   -L\,\ell(\overline{X})\,\textbf{T}(\overline{a},\overline{Y}, \overline{Z})+L^{2}\,\textbf{T}(
   T(\overline{a},\overline{X}),\overline{Y}, \overline{Z})\} \nonumber\\
  && +f_{2}\{L\,\ell(\overline{X})\,\textbf{T}(\overline{b},\overline{Y}, \overline{Z})-L^{2}\,\textbf{T}(
   T(\overline{b},\overline{X}),\overline{Y}, \overline{Z})\}  +u(\overline{X})\textbf{T}(\varphi_{2}(\overline{\eta}),\overline{Y}, \overline{Z}) \nonumber\\
  && +L\,\ell(\varphi_{1}(\overline{X}))\,\textbf{T}(\overline{u},\overline{Y}, \overline{Z})-L\,\ell(\varphi_{1}(\overline{\eta}))\, \textbf{T}(T(\overline{u},\overline{X}),\overline{Y}, \overline{Z})\nonumber\\
  &&+u(\overline{\eta})\,  \textbf{T}(T((\varphi_{1}-\varphi_{2})(\overline{\eta}) , \overline{X}) -\varphi_{1}(\overline{X}),\overline{Y}, \overline{Z}) .
  \end{eqnarray}
In addition, the Berwald connection $D^{\circ}$ can be written as follows
  \begin{eqnarray}\label{2.11}
 {D^{\circ}} _{\gamma \overline{X}_{t}}\overline{Y}&=&f_{1}\,\{A(\overline{X})\,{D^{\circ}}_{\gamma \overline{\eta}}\,\overline{Y}+
   A(\overline{\eta})\,{D^{\circ}} _{\gamma \overline{X}}\overline{Y}-L\,\ell(\overline{X}){D^{\circ}} _{\gamma \overline{a}}\overline{Y}
   +L^{2}\,{D^{\circ}} _{\gamma T(\overline{a},\overline{X})}\overline{Y}\} \nonumber\\
  &&+ f_{2}\{L\,\ell(\overline{X})\,{D^{\circ}} _{\gamma \overline{b}}\overline{Y}
   -L^{2}\,{D^{\circ}} _{\gamma T(\overline{b},\overline{X})}\overline{Y}\} -L\,\ell(\varphi_{1}(\overline{\eta}))\,{D^{\circ}} _{\gamma T(\overline{u},\overline{X})}\overline{Y}
\nonumber\\
&&   +u(\overline{X})\,{D^{\circ}}_{\gamma\varphi_{2}(\overline{\eta})}\overline{Y} -u(\overline{\eta})\,{D^{\circ}} _{\gamma\varphi_{1}(\overline{X}) -{\gamma T((\varphi_{1}-\varphi_{2})(\overline{\eta}), \overline{X})}}\overline{Y}\nonumber \\&&  +L\,\ell(\varphi_{1}(\overline{X}))\,{D^{\circ}} _{\gamma\overline{u}}\overline{Y}.
\end{eqnarray}
Using the above three relations  (\ref{B1}) - (\ref{2.11}), together with the formula $\nabla_{\gamma\, \overline{X}} \overline{Y} = {D^{\circ}}_{\gamma\, \overline{X}} \overline{Y}+T(\overline{X},\overline{Y})$ \cite{r92},  we conclude from  (\ref{eq.t4})  that  {\fontsize{11.3}{11.3}\selectfont
 \begin{eqnarray}\label{22}
\overline{D} _{{\bar{\beta}\,} \overline{X}} \overline{Y}&=&\nabla _{{\beta}\, \overline{X}} \overline{Y}
+f_{1}\{A(\overline{\eta})\,\nabla_{\gamma\, \overline{X}} \overline{Y} +A(\overline{X})\,\nabla_{\gamma\, \overline{\eta}} \overline{Y}
-A(\overline{X})\,\overline{Y}-A(\overline{Y})\,\overline{X}+\overline{a}\,g(\overline{X},\overline{Y})\nonumber\\
&&-L\,\{\ell(\overline{X})\,{\nabla}_{\gamma \overline{a}}\overline{Y} +\ell(\overline{Y})\,T(\overline{a},\overline{X})\}
   +L^{2}\{ \nabla_{\gamma T(\overline{a},\overline{X})}\overline{Y} +S(\overline{X},\overline{a})\overline{Y}\}+\textbf{T}(\overline{a},\overline{X},\overline{Y})\overline{\eta}\}\nonumber \\
&&-f_{2}\,\{g(\overline{X},\overline{Y})\overline{b}-L\,\{\ell(\overline{X})\,{\nabla}_{\gamma \overline{b}}\overline{Y} +\ell(\overline{Y})\,T(\overline{b},\overline{X})\}
   +L^{2}\,\{S(\overline{X},\overline{b})\overline{Y}+ \nabla_{\gamma T(\overline{b},\overline{X})}\overline{Y}\}\nonumber\\
   &&+
\textbf{T}(\overline{b},\overline{X},\overline{Y})\overline{\eta}\}  -g(\varphi_{1}(\overline{X}),\overline{Y})\,\overline{u}+u(\overline{X})\{ \nabla_{\gamma\varphi_{2}(\overline{\eta})}\overline{Y} -\varphi_{2}(\overline{Y})\}\nonumber\\
&&+u(\overline{Y})\{\varphi_{1}(\overline{X})+T(\varphi_{2}(\overline{\eta}),\overline{X})\}
+ L\,\ell(\varphi_{1}(\overline{X}))\, \nabla_{\gamma\overline{u}}\overline{Y}+ L\,\ell(\varphi_{1}(\overline{Y}))\,T(\overline{u},\overline{X})
\nonumber\\
&& +L\,\ell(\varphi_{1}(\overline{\eta}))\,\{S(\overline{u},\overline{X})\overline{Y}-\nabla_{\gamma T(\overline{u},\overline{X})}\overline{Y} \} -\mathbf{T}(\overline{u},\overline{X},\overline{Y})\,\varphi_{1}(\overline{\eta})-
\mathbf{T}(\varphi_{2}(\overline{\eta}),\overline{X},\overline{Y})\,\overline{u}\nonumber\\
&&
+ u(\overline{\eta})\{\varphi_{1}(T(\overline{X},\overline{Y}))-T(\varphi_{1}(\overline{Y}),\overline{X})
+S(\overline{X},(\varphi_{1}- \varphi_{2})(\overline{\eta}))\overline{Y}
-\nabla_{\gamma\varphi_{1}(\overline{X})}\overline{Y}\nonumber\\
&& +\nabla_{\gamma T((\varphi_{1}-\varphi_{2})(\overline{\eta}), \overline{X})}\overline{Y}\}.
\end{eqnarray}}
\par Consequently, from (\ref{11}) and (\ref{22}),  taking into account (\ref{B1}), the full
expression of $\overline{D}_{X}\overline{Y}$ in terms of Cartan connection is the following
\begin{eqnarray}\label{33}
\overline{D}_{X} \overline{Y}&=&\nabla _{X} \overline{Y}
+f_{1}\,\{g(\rho {X},\overline{Y})\,\overline{a}-A(\rho {X})\,\overline{Y}-A(\overline{Y})\,\rho {X}-L\,\ell(\overline{Y})\,T(\overline{a},\rho {X})\nonumber\\
&&+\textbf{T}(\overline{a},\rho {X},\overline{Y})\,\overline{\eta}
+L^{2}\,S(\rho{X},\overline{a})\overline{Y}\}-\textbf{T}(\overline{u},\rho {X},\overline{Y})\,\varphi_{1}(\overline{\eta})-\mathbf{T}(\varphi_{2}(\overline{\eta}),\rho {X},\overline{Y})\,\overline{u}\nonumber \\
&&-f_{2}\,\{g(\rho {X},\overline{Y})\,\overline{b}-L\,\ell(\overline{Y})\,T(\overline{b},\rho {X})
+\textbf{T}(\overline{b},\rho {X},\overline{Y})\,\overline{\eta} +L^{2}S(\rho{X},\overline{b})\overline{Y}\}\nonumber \\
&&-g(\varphi_{1}(\rho {X}),\overline{Y})\,\overline{u}-u(\rho {X})\,\varphi_{2}(\overline{Y})+u(\overline{Y})\,\{\varphi_{1}(\rho {X}) +T(\varphi_{2}(\overline{\eta}),\rho {X})\}\nonumber\\
&&
-u(\overline{\eta})\,\{T(\varphi_{1}(\overline{Y})  ,\rho {X})+S((\varphi_{1}-\varphi_{2})(\overline{\eta}),\rho {X})\overline{Y}
+\varphi_{1}(T(\rho {X},\overline{Y}))\} \nonumber
\\
&&
+ L\,\ell(\varphi_{1}(\overline{Y}))T(\overline{u},\rho {X})
 +L\,\ell(\varphi_{1}(\overline{\eta}))\,
S(\overline{u},\rho {X})\overline{Y}.
\end{eqnarray}
Hence, ${\overline{D}}_{X}\overline{Y}$ is uniquely determined
by the right-hand side of (\ref{33}).

  \vspace{7pt}
  \par
 In order to prove the existence of $\overline{D}$, just define $\overline{D}$ by the above formula. Then, it is easy to check that  $\overline{D}$ is a regular Finsler connection that satisfies the conditions (I) - (IV). This completes the proof.
\end{proof}

\begin{rem}\label{remm1} It is worth mentioning that  the  connection  $GC\Gamma$ is  the Finslerian version  of the Riemannian Tripathi connection  \emph{\cite{New}} and generalizes the local study provided in \emph{\cite{amr20}}.
\end{rem}

\begin{cor}\label{th.v3} The  $GC\Gamma$-connection  $\overline{D}$ and  the Cartan connection $\nabla$ are related by
$$\overline{D} _{X} \overline{Y}=\nabla _{X} \overline{Y}+ N(\rho X, \overline{Y}),$$ where
{\fontsize{11.48}{11.48}\selectfont
  \begin{eqnarray*}\label{B3a}
N(\rho X, \overline{Y})&=&f_{1} \,\{g(\rho {X},\overline{Y})\,\overline{a}-A(\rho {X})\overline{Y}-A(\overline{Y})\,\rho {X}-L\,\ell(\overline{Y})\,T(\overline{a},\rho {X})+\textbf{T}(\overline{a},\rho {X},\overline{Y})\overline{\eta}
\\
&&+L^{2}\,S(\rho{X},\overline{a})\overline{Y} \} -f_{2}\{g(\rho {X},\overline{Y})\,\overline{b}-L\,\ell(\overline{Y})\,T(\overline{b},\rho {X})
+\textbf{T}(\overline{b},\rho {X},\overline{Y})\,\overline{\eta} \\
&&+L^{2}\,S(\rho{X},\overline{b}\overline{Y}) \}-\{g(\varphi_{1}(\rho {X}),\overline{Y})+
\mathbf{T}(\varphi_{2}(\overline{\eta}),\rho {X},\overline{Y})\}\overline{u}-u(\rho {X})\,\varphi_{2}(\overline{Y})
\\
&&+ L\,\ell(\varphi_{1}(\overline{Y}))T(\overline{u},\rho {X})-u(\overline{\eta})\{S((\varphi_{1}-\varphi_{2} )(\overline{\eta}),\rho {X} ) \overline{Y} +T(\varphi_{1}(\overline{Y}),\rho {X})\\
&&
  -\varphi_{1}(T(\rho {X},\overline{Y}))\} +u(\overline{Y})\{T(\varphi_{2}(\overline{\eta}),\overline{X})+\varphi_{1}(\rho {X})\} \\
&&+L\,\ell(\varphi_{1}(\overline{\eta}))\,S(\overline{u},\rho {X})\overline{Y}
-\textbf{T}(\overline{u},\rho {X},\overline{Y})\,\varphi_{1}(\overline{\eta}) .
\end{eqnarray*}}
\end{cor}
\begin{prop} \label{th.v2} Consider a Finsler manifold $(M,L)$ and let $\overline{D}$ be the Tripathi Finsler connection.   Then, \vspace{-0.1cm}
\begin{description}
  \item[(a)]The  canonical spray $\bar{G}$ associated with $\overline{D}$ is given by  \vspace{-0.3cm}
$$\bar{G}= G+f_{1}\,\{2A(\overline{\eta})\,\gamma\overline{\eta}-L^{2}\,\gamma\overline{a}\}+f_{2}\,L^{2}\,\gamma\overline{b} +L\, \ell(\varphi_{1}(\overline{\eta}))\,\gamma\overline{u}-u(\overline{\eta})\,\gamma(\varphi_{1}-\varphi_{2})(\overline{\eta}),\vspace{-0.3cm}$$
where $G$ is the geodesic spray of the Finsler structure.
\item[(b)]The canonical nonlinear connection $ \overline{\Gamma}$ associated with $\overline{D}$  is given by:\vspace{-0.3cm}
{\fontsize{11.2}{11.2}\selectfont
\begin{eqnarray*}
 \overline{\Gamma}(X):=\Gamma(X)&+&2\Big( f_{1}\{A(\rho{X})\,\gamma\overline{\eta}
   +A(\overline{\eta})\,J{X}-L\,\ell(\rho{X})\,\gamma \overline{a} +L^{2}\gamma T(\overline{a},\rho{X})\} \nonumber\\
   &&  + f_{2}\,\{L\,\ell(\rho{X})\,\gamma \overline{b} -L^{2}\gamma T(\overline{b},\rho{X})\} +L\,\ell(\varphi_{1}(\rho{X}))\,\gamma\overline{u}-u(\overline{\eta})\gamma\varphi_{1}(\rho{X})\nonumber\\
  && +u(\rho{X}) \gamma \varphi_{2}(\overline{\eta})- L\ell(\varphi_{1}(\overline{\eta}))\gamma T(\overline{u},\rho{X})+u(\overline{\eta})\,
 \gamma T((\varphi_{1}- \varphi_{2})(\overline{\eta}), \rho{X})\Big), \end{eqnarray*}}
 where $\Gamma$ is the Barthel connection of $(M,L)$. \end{description}
     \end{prop}

\begin{proof}
$\textbf{(a)}$ It follows from \eqref{diff}, by replacing $\overline{X}$ by $\overline{\eta}$, and using \eqref{B}.
\\
 $\,\textbf{(b)}$ The expression of $\overline{\Gamma}$ is obtained by  applying Lemma \ref{eqv.}, taking into account Equations \eqref{diff} and \eqref{B1}.
\end{proof}

\begin{prop}\label{prop.v1} Let $\overline{D}$ be the Tripathi Finsler connection.  Then, the following hold:
\begin{description}
  \item[(a)] the $(h)hv$-torsion $\overline{T}$ of $\overline{D}$ coincides with the $(h)hv$-torsion ${T}$ of Cartan connection.\vspace{-0.3cm}
  \item[(b)] the $(h)h$-torsion ${\overline{Q}}$ of $\overline{D}$ has the form:
  $\overline{Q}(\overline{X},\overline{Y})= u(\overline{Y}) \varphi(\overline{X})-u(\overline{X}) \varphi(\overline{Y}).$

  \item[(c)]the $(v)v$-torsion $\widehat{\overline{S}}$ of $\overline{D}$ vanishes identically.

  \item[(d)]the $(v)hv$-torsion $\widehat{\overline{P}}$ of $\overline{D}$ has the form:
  $$\widehat{\overline{P}}(\overline{X},\overline{Y})= \widehat{P}(\overline{X},\overline{Y}) -\nabla_{{\gamma} \overline{Y}} \overline{X}_{t} -N(\overline{X},\overline{Y})+
   N(\rho [{\beta} \overline{X}, \gamma \overline{Y}], \overline{\eta}).$$\vspace{-0.9cm}
  \item[(e)]the $(v)h$-torsion $\widehat{\overline{R}}$ of $\overline{D}$ has the form:
  $$\widehat{\overline{R}}(\overline{X},\overline{Y})= \widehat{{R}}(\overline{X},\overline{Y})+ N(\rho [{\bar{\beta}}\, \overline{X}, {\bar{\beta}}\, \overline{Y}], \overline{\eta})+K([{\beta} \overline{X}, \gamma \overline{Y}_{t}]+[\gamma \overline{X}_{t}, \beta \overline{Y}]+[{\gamma} \overline{X}_{t}, \gamma \overline{Y}_{t}]),$$
\end{description}\vspace{-0.3cm}
where $\widehat{{P}}$ and $\widehat{{R}}$ are the $(v)hv$ and $(v)h$ torsions of the Cartan connection, respectively. 

\end{prop}
\begin{proof}\textbf{\textbf{(a)}} Follows from the definition of $\overline{T}$, together with  Equations  (\ref{eq.t2}) and   (\ref{11}).

\vspace{4pt}
 \noindent\textbf{\textbf{(b)}} Follows directly by condition \textbf{(III)} of Theorem \ref{th.v1}.

 \vspace{4pt}
 \noindent\textbf{\textbf{(c)}} Using \textbf{(a) }above and   \cite[Proposition 2.5]{r96}, we obtain
 \begin{eqnarray}\label{Scurvature}
   \overline{S}(\overline{X},\overline{Y})\overline{Z} &=& (\overline{D}_{\gamma \overline{Y}}\overline{T})(\overline{X},\overline{Z})-
 (\overline{D}_{\gamma \overline{X}}\overline{T})(\overline{Y},\overline{Z})  +\overline{T}(\overline{X},\overline{T}(\overline{Y},\overline{Z}))\nonumber\\ 
 && -\overline{T}(\overline{Y},\overline{T}(\overline{X},\overline{Z})) +\overline{T}(\widehat{\overline{S}}(\overline{Y},\overline{X}),\overline{Z}).
 \end{eqnarray}
Setting $\overline{Z}=\overline{\eta}$ into (\ref{Scurvature}), taking into account  \textbf{(a)} together with the properties of $T$ and the fact that  $ \bar{K}\circ\gamma=id_{\cp}$, the result follows.

  \vspace{4pt}
 \noindent\textbf{\textbf{(d)}} According to Corollary  \ref{th.v3} and   Proposition \ref{th.v2} together with $ \bar{K}\circ\bar{\beta}=0$, we get
 \begin{eqnarray*}
   \widehat{\overline{P}}(\overline{X},\overline{Y}) &=& -\overline{D}_{\bar{\beta}\, \overline{X}}  \overline{D}_{\gamma \overline{Y}}\, \overline{\eta}
   +  \overline{D}_{\gamma \overline{Y}}\overline{D}_{\bar{\beta}\, \overline{X}} \, \overline{\eta}+ \overline{D}_{[\bar{\beta}\, \overline{X}, \gamma \overline{Y}]} \,\overline{\eta} 
   =-\overline{D}_{\bar{\beta}\, \overline{X}} \overline{Y}+ \overline{D}_{[\bar{\beta}\, \overline{X}, \gamma \overline{Y}]} \overline{\eta}\\
   &=& -\nabla_{{\beta} \,\overline{X}} \overline{Y}+ \nabla_{[{\beta}\, \overline{X}, \gamma \overline{Y}]} \overline{\eta}
    -\nabla_{{\gamma} \overline{X}_{t}} \overline{Y}+ \nabla_{[{\gamma} \overline{X}_{t}, \gamma \overline{Y}]} \overline{\eta}-N(\overline{X},\overline{Y})+    N(\rho [{\bar{\beta}}\, \overline{X}, \gamma \overline{Y}], \overline{\eta})\\
   &=& \widehat{P}(\overline{X},\overline{Y}) -\nabla_{{\gamma} \overline{Y}} \overline{X}_{t} -N(\overline{X},\overline{Y})+
   N(\rho [{\beta}\, \overline{X}, \gamma \overline{Y}], \overline{\eta}).
  \end{eqnarray*}
 Hence, the result follows by taking into account  \textbf{(c)} above.

 \vspace{4pt}
 \noindent\textbf{\textbf{(e)}} The proof is similar to that of \textbf{(d)} above.
 \end{proof}


\begin{prop}\label{7.pp.7} Let $\overline{S}$,  $\overline{P}$ and $\overline{R}$ \emph{(}$S, P$ and $R$\emph{)} be the v-, hv- and h-curvatures of the $GC\Gamma$-connection  $\overline{D}$ \emph{(}the Cartan connection $\nabla$\emph{)},
then we have~ \vspace{-0.1cm}
\begin{description}

  \item[(a)]${{\overline{S}}}(\overline{X},\overline{Y})\overline{Z}= S(\overline{X},\overline{Y})\overline{Z}$,

  \item[(b)]${{\overline{P}}}(\overline{X},\overline{Y})\overline{Z}= P(\overline{X},\overline{Y})\overline{Z}
  +  (\nabla_{\gamma \overline{Y}}N)(\overline{X},\overline{Z})+N({T}(\overline{Y},\overline{X}),\overline{Z})$\\
  ${\qquad\qquad\ \ \ }+f_{1}\{A(\overline{\eta})\,S(\overline{X},\overline{Y})\overline{Z}-L\,\ell(\overline{X})\,
  S(\overline{a},\overline{Y})\overline{Z}+L^{2}\,S(T(\overline{a},\overline{X}),\overline{Y})\overline{Z} \}$\\
  ${\qquad\qquad\ \ \ }+f_{2}\{L\,\ell(\overline{X})\, S(\overline{b},\overline{Y})\overline{Z}+L^{2}\,S(T(\overline{b},\overline{X}),\overline{Y})\overline{Z} \}+L\,\ell(\varphi_{1}(\overline{X}))\,S(\overline{u},\overline{Y})\overline{Z}$\\
  ${\qquad\qquad\ \ \ }+u(\overline{X})\,S(\varphi_{2}(\overline{\eta}),\overline{Y})\overline{Z}
 -L\,\ell(\varphi_{1}(\overline{\eta}))\, S(T(\overline{u},\overline{X}),\overline{Y})\overline{Z}  $\\
  ${\qquad\qquad\ \ \ }+
  u(\overline{\eta})\,\{S(T((\varphi_{1}-\varphi_{2})(\overline{\eta}),\overline{X}),\overline{Y})\overline{Z}-S(\varphi_{1}(\overline{X}),\overline{Y})\overline{Z} \}$,

  \item[(c)]${{\overline{R}}}(\overline{X},\overline{Y})\overline{Z}= R(\overline{X},\overline{Y})\overline{Z}+
  P(\overline{X},\overline{Y}_{t})\overline{Z}-P(\overline{Y},\overline{X}_{t})\overline{Z}+S(\overline{X}_{t},\overline{Y}_{t})\overline{Z}$\\
      ${\qquad\qquad\ \ \ }
   +\mathfrak{U}_{\overline{X}, \overline{Y}}\{(\nabla_{\bar{\beta}\,\overline{Y}}N)(\overline{X},\overline{Z})
   +N(\overline{Y}, N(\overline{X}, \overline{Z}))+N(T(\overline{Y}_{t},\overline{X}),\overline{Z})\}$,
\end{description}
where $\mathfrak{U}_{X, Y}\{B(X,Y)\}:=B(X,Y)-B(Y,X)$.
\end{prop}
\begin{proof}\textbf{(a)} Follows from \eqref{11}. \textbf{(b)} and \textbf{(c)} follow from both \eqref{11} and \eqref{22}.
\end{proof}
\begin{prop}\label{pp.1} For the $GC\Gamma$-connection, the following identities~hold{\em:}\,
\begin{description}

\item[\textbf{(a)}]$\overline{P}(\overline{X},\overline{Y})\overline{Z}-\overline{P}(\overline{Z},\overline{Y})\overline{X}=
(\overline{D}_{\bar{\beta} \overline{Z}}{T})(\overline{Y},\overline{X})-(\overline{D}_{\bar{\beta}
\overline{X}}{T})(\overline{Y},\overline{Z})-(\overline{D}_{\gamma
\overline{Y}}\overline{Q})(\overline{Z},\overline{X})$\\
$ { \qquad\qquad\qquad \qquad\qquad\ \ }+{T}(\overline{Y},
\overline{Q}(\overline{Z},\overline{X}))-{T}(\widehat{\overline{P}}(\overline{Z},\overline{Y}),\overline{X})+
{T}(\widehat{\overline{P}}(\overline{X},\overline{Y}),\overline{Z})$
\\ $ { \qquad\qquad\qquad \qquad\qquad\ \ }-\overline{Q}(\overline{Z},
{T}(\overline{Y},\overline{X})) +\overline{Q}(\overline{X},
{T}(\overline{Y},\overline{Z})),$

 \item[\textbf{(b)}]$\mathfrak{S}_{\overline{X},\overline{Y},\overline{Z}}
\{\overline{R}(\overline{X},
\overline{Y})\overline{Z}-T(\widehat{\overline{R}}(\overline{X},\overline{Y}),\overline{Z})\}=
\mathfrak{S}_{\overline{X},\overline{Y},\overline{Z}}
\{\overline{Q}(\overline{X},\overline{Q}(\overline{Y},\overline{Z}))-(\overline{D}_{\bar{\beta}\,
\overline{X}}\overline{Q})(\overline{Y},\overline{Z})\}$,
\item[\textbf{(c)}]
$(\overline{D}_{\bar{\beta}\,\overline{Z}}S)(\overline{X},\overline{Y},\overline{W})
-\overline{P}(\overline{Z},\widehat{S}(\overline{X},\overline{Y}))\overline{W}=$\\$
=\mathfrak{U}_{\overline{X}, \overline{Y}}\{
(\overline{D}_{\gamma \overline{Y}}\overline{P})(\overline{Z},\overline{X},
\overline{W})+\overline{P}(T(\overline{X},\overline{Z}),\overline{Y})\overline{W}+S(\widehat{\overline{P}}(\overline{Z},\overline{X}),\overline{Y})\overline{W}\},$
\item[\textbf{(d)}]
$(\overline{D}_{\gamma\overline{X}}\overline{R})(\overline{Y},\overline{Z},\overline{W})=
S(\widehat{\overline{R}}(\overline{Y},\overline{Z}),\overline{X})\overline{W}- \overline{P}(\overline{Q}(\overline{Y},\overline{Z}),\overline{X})\overline{W}
+$\\$+
\mathfrak{U}_{\overline{Z}, \overline{Y}}\{ (\overline{D}_{\bar{\beta}   \overline{Z}}\overline{P})(\overline{Y},\overline{X},\overline{W})   + \overline{P}(\overline{Z},\widehat{\overline{P}}(\overline{Y},\overline{X}))\overline{W}+\overline{R}(T(\overline{X},\overline{Z}),\overline{Y})\overline{W}\},$
\item[\textbf{(e)}]
$\mathfrak{S}_{\overline{X},\overline{Y},\overline{Z}} \{(\overline{D}_{\bar{\beta}
\overline{X}}\overline{R})(\overline{Y},
\overline{Z},\overline{W})+\overline{P}(\overline{X},\widehat{\overline{R}}(\overline{Y},\overline{Z}))\overline{W}
+\overline{R}(\overline{Q}(\overline{X},\overline{Y}),\overline{Z})\overline{W}\}=0$.
\end{description}
\end{prop}
\begin{proof} It results from  \cite[Propositions 2.5, 2.6]{r96}, taking into account  Corollary \ref{th.v3} and  Propositions \ref{prop.v1}, \ref{7.pp.7} above.
\end{proof}
\vspace{-0.5 cm}
\section{A generalization of the four celebrated Finsler connections}
 This section is devoted to constructing new  Finsler connections from our general connection $G C{\Gamma}$ by means of the ${P}^{1}$-process and ${C}$-process introduced by Matsumoto \cite{Ma.book2}. First,  let us denote Cartan, Berwald, Hasiguashi and Rund-Chern connections by $C\Gamma,\,\, B\Gamma,\,\,H\Gamma$ and $R\Gamma$,  respectively.
\begin{defn} \label{def.G.}  Let  $\overline{D}$ be the $GC\Gamma$-connection.  The process of adding  the associated $(v)$hv-torsion $\overline{P}(\overline{Y},\overline{X})$ to  the horizontal part $\overline{D}_{\bar{\beta}\, \overline{X}}\overline{Y}$  of  $GC\Gamma$ is called the $\overline{P}^{1}$-process. Moreover, the process of subtracting  the associated $(h)$hv-torsion $\overline{T}(\overline{X},\overline{Y})$ from  the vertical part $\overline{D}_{\gamma \overline{X}}\overline{Y}$  of  $GC\Gamma$ is called $\overline{C}$-process.
\end{defn}

  \begin{thm}\label{G_F_C}
By means of the  $\overline{P}^{1}\!$-process and $\overline{C}$-process, we have:
\begin{description}
  \item[(a)] The  $ \overline{P}^{1}$-process of $GC\Gamma$ yields  a general  Hashiguchi  connection \emph{($GH \Gamma$)}.
    \item[(b)] The $ \overline{C}$-process of $GC\Gamma$ yields   a general Rund-Chern connection \emph{($GR \Gamma$)}.
   \item[(c)] The $ \overline{P}^{1}$-process followed by the
$\overline{C}$-process  of $GC\Gamma$ \emph{(}or the $\overline{C}$-process followed by the $ \overline{P}^{1}$-process\emph{) }yields  a general  Berwald  connection \emph{($GB \Gamma$)}.
\end{description}
  \end{thm}

Now, we define a vanishing condition (or in short VC) by setting $f_{1}=f_{2}=u=0$.  Based on Theorem \ref{G_F_C} and Corollary \ref{th.v3}, if the VC is  satisfied, we obtain the following diagram:
{\large{{\[
\begin{CD}
       GB \Gamma @< \overline{C}\text{-proc.} << GH \Gamma @< \overline{P}^{1}\text{-proc.} <<  G C\Gamma @> \overline{C}\text{-proc.}  >> GR\Gamma @> \overline{P}^{1}\text{-proc.}  >> GB\Gamma\\
    @V \text{VC}  VV @V \text{VC} VV   @V \text{VC}  VV @VV \text{VC}  V @VV \text{VC}  V\\
        B \Gamma @< {C}\text{-proc.} <<  H \Gamma @< {P}^{1}\text{-proc.} <<  C\Gamma @> {C}\text{-proc.}  >> R\Gamma @> {P}^{1}\text{-proc.}  >> B\Gamma
\end{CD}
\]}}}
The arrows of the second row arise from the usual $P^{1}$-process and $C$-process and they are well known \cite{Ma.book2}. The arrows of the first row arise from Theorem \ref{G_F_C} and they are completely new. Moreover, the connections of the second row come from the general connections of the first row under the VC. 

\section{Special cases}

In the present section, we give some important  particular cases (26 cases) of  our connection $\overline{D}$  which result from  certain choices of $f_{1}, f_{2}, A, B,u,\varphi_{1}, \varphi_{2}$. Some of the following cases have been already studied in the context of Finsler geometry while  many others have not.

\begin{flushleft}
\textbf{Generalized quarter-symmetric recurrent metric Finsler connection: $ A=B,\,\,f_{1}=1-t ,\, f_{2}= -t, \,\, t \in \mathbb{R}$ }
\end{flushleft}
\begin{description}
  \item{(1)} In this case,  $\overline{D}$  reduces to  a Finslerian version of the connection introduced in \cite{gen quart} and is written in the form:
 \begin{eqnarray*}
\overline{D} _{X} \overline{Y}&=&\nabla _{X} \overline{Y} -(1-t) \,\{A(\rho {X})\overline{Y}+A(\overline{Y})\,\rho {X} \} + g(\rho {X},\overline{Y})\,\overline{a}-L\,\ell(\overline{Y})\,T(\overline{a},\rho {X})
\nonumber\\
&&+L^{2}\,S(\rho{X},\overline{a})\overline{Y} -\{g(\varphi_{1}(\rho {X}),\overline{Y})+
\mathbf{T}(\varphi_{2}(\overline{\eta}),\rho {X},\overline{Y})\}\overline{u}-u(\rho {X})\,\varphi_{2}(\overline{Y})
\nonumber\\
&&+ L\,\ell(\varphi_{1}(\overline{Y}))T(\overline{u},\rho {X})-u(\overline{\eta})\{S((\varphi_{1}-\varphi_{2} )(\overline{\eta}),\rho {X} ) \overline{Y} +T(\varphi_{1}(\overline{Y}),\rho {X})\nonumber\\
&&
  -\varphi_{1}(T(\rho {X},\overline{Y}))\} +u(\overline{Y})\{T(\varphi_{2}(\overline{\eta}),\overline{X})+\varphi_{1}(\rho {X})\} +\textbf{T}(\overline{a},\rho {X},\overline{Y})\overline{\eta}\\
&&+L\,\ell(\varphi_{1}(\overline{\eta}))\,S(\overline{u},\rho {X})\overline{Y}
-\textbf{T}(\overline{u},\rho {X},\overline{Y})\,\varphi_{1}(\overline{\eta}) .
\end{eqnarray*}
\end{description}
\textbf{Quarter-symmetric metric Finsler connection: \,$f_{1}=f_{2}=0$}

\begin{description}
  \item{(2)} $\overline{D}$ becomes a Finslerian version of the Riemannian connection given in \cite[formula (3.3)]{ Yano.1}, that is,
\begin{eqnarray*}
\overline{D} _{X} \overline{Y}&=&\nabla _{X} \overline{Y} -\{g(\varphi_{1}(\rho {X}),\overline{Y})+
\mathbf{T}(\varphi_{2}(\overline{\eta}),\rho {X},\overline{Y})\}\overline{u}-u(\rho {X})\,\varphi_{2}(\overline{Y})
\nonumber\\
&&+ L\,\ell(\varphi_{1}(\overline{Y}))T(\overline{u},\rho {X})-u(\overline{\eta})\{S((\varphi_{1}-\varphi_{2} )(\overline{\eta}),\rho {X} ) \overline{Y} +T(\varphi_{1}(\overline{Y}),\rho {X})\nonumber\\
&&
  -\varphi_{1}(T(\rho {X},\overline{Y}))\} +u(\overline{Y})\{T(\varphi_{2}(\overline{\eta}),\overline{X})+\varphi_{1}(\rho {X})\} \nonumber\\
&&+L\,\ell(\varphi_{1}(\overline{\eta}))\,S(\overline{u},\rho {X})\overline{Y}
-\textbf{T}(\overline{u},\rho {X},\overline{Y})\,\varphi_{1}(\overline{\eta})
\nonumber .
\end{eqnarray*}
\item{(3)} In addition, when  $\varphi=\emph{Ric}_{o}$, then we obtain a Finslerian version of the Ricci quarter-symmetric metric connection appeared in \cite[formula (2.2)]{Pandey.1}.
\item{(4)} When $\varphi_{2}=0$, we get a Finslerian version of the quarter-symmetric metric connection presented in \cite[formula (1.6)]{Pandey.1}. Thus,   $\overline{D}$ has the form:
\begin{eqnarray*}
\overline{D} _{X} \overline{Y}&=&\nabla _{X} \overline{Y} -g(\varphi_{1}(\rho {X}),\overline{Y})\,\overline{u}+ L\,\ell(\varphi_{1}(\overline{Y}))T(\overline{u},\rho {X})-u(\overline{\eta})\{S(\varphi_{1}(\overline{\eta}),\rho {X} ) \overline{Y}\nonumber\\
&& +T(\varphi_{1}(\overline{Y}),\rho {X})
  -\varphi_{1}(T(\rho {X},\overline{Y}))\} +u(\overline{Y})\,\varphi_{1}(\rho {X})\nonumber\\
&& +L\,\ell(\varphi_{1}(\overline{\eta}))\,S(\overline{u},\rho {X})\overline{Y}
-\textbf{T}(\overline{u},\rho {X},\overline{Y})\,\varphi_{1}(\overline{\eta})
\nonumber .
\end{eqnarray*}
\item{(5)} When $\varphi_{1}=0,\,\,\overline{D}$ reduces to a Finslerian version of the connection  given in \cite[formula (3.6)]{Yano.1}. That is,
$$\begin{aligned}
\overline{D} _{X} \overline{Y}&=\nabla _{X} \overline{Y} -\mathbf{T}(\varphi_{2}(\overline{\eta}),\rho {X},\overline{Y})\,\overline{u}-u(\rho {X})\varphi_{2}(\overline{Y})
+u(\overline{Y})\,T(\varphi_{2}(\overline{\eta}),\overline{X})\nonumber\\
&+u(\overline{\eta})\,S(\varphi_{2}(\overline{\eta}),\rho {X})\overline{Y}
\nonumber .
\end{aligned}$$
\end{description}

\newpage
\begin{flushleft}
\textbf{Quarter-symmetric non-metric Finsler connection: $f_{1}\neq 0,\,f_{2}=0$}
\end{flushleft}

\begin{description}
  \item{(6)} If $f_{1}=\frac{1}{2}$,  then we obtain an intrinsic formula of the quarter-symmetric $h$-recurrent Finsler connection presented in \cite{amr20}.  Thereby, $\overline{D}$ becomes \vspace{-0.25 cm}
\begin{eqnarray*}
\overline{D} _{X} \overline{Y}&=&\nabla _{X} \overline{Y}+\frac{1}{2}\,\{g(\rho {X},\overline{Y})\,\overline{a}-A(\rho {X})\overline{Y}-A(\overline{Y})\,\rho {X}-L\,\ell(\overline{Y})\,T(\overline{a},\rho {X})
\nonumber\\
&&+L^{2}\,S(\rho{X},\overline{a})\overline{Y} \} -\{g(\varphi_{1}(\rho {X}),\overline{Y})+
\mathbf{T}(\varphi_{2}(\overline{\eta}),\rho {X},\overline{Y})\}\overline{u}\nonumber\\
&&-u(\rho {X})\,\varphi_{2}(\overline{Y})
+ L\,\ell(\varphi_{1}(\overline{Y}))T(\overline{u},\rho {X})-u(\overline{\eta})\{S((\varphi_{1}-\varphi_{2} )(\overline{\eta}),\rho {X} ) \overline{Y} \nonumber\\
&&+T(\varphi_{1}(\overline{Y}),\rho {X})
  -\varphi_{1}(T(\rho {X},\overline{Y}))\} +u(\overline{Y})\{T(\varphi_{2}(\overline{\eta}),\overline{X})+\varphi_{1}(\rho {X})\} \nonumber\\
&&+L\,\ell(\varphi_{1}(\overline{\eta}))\,S(\overline{u},\rho {X})\overline{Y}
-\textbf{T}(\overline{u},\rho {X},\overline{Y})\,\varphi_{1}(\overline{\eta})+\frac{1}{2}\textbf{T}(\overline{a},\rho {X},\overline{Y})\overline{\eta} \nonumber .
\end{eqnarray*}

 \item{(7)} When $\varphi_{2}=0,\,\,\overline{D}$  reduces to a Finslerian version of the connection given in \cite[\S 4.2 (4)]{New}. Then,   $\overline{D}$ can be written as follows
\begin{eqnarray*}
\overline{D} _{X} \overline{Y}&=&\nabla _{X} \overline{Y}+f_{1} \,\{g(\rho {X},\overline{Y})\,\overline{a}-A(\rho {X})\overline{Y}-A(\overline{Y})\,\rho {X}-L\,\ell(\overline{Y})\,T(\overline{a},\rho {X})\nonumber\\
&&+\textbf{T}(\overline{a},\rho {X},\overline{Y})\overline{\eta}
+L^{2}\,S(\rho{X},\overline{a})\overline{Y} \} + L\,\ell(\varphi_{1}(\overline{Y}))T(\overline{u},\rho {X})+u(\overline{Y})\varphi_{1}(\rho {X})
\nonumber\\
&&
-u(\overline{\eta})\{S(\varphi_{1}(\overline{\eta}),\rho {X} ) \overline{Y} +T(\varphi_{1}(\overline{Y}),\rho {X}) -\varphi_{1}(T(\rho {X},\overline{Y}))\}
 \nonumber\\
&&+L\,\ell(\varphi_{1}(\overline{\eta}))\,S(\overline{u},\rho {X})\overline{Y}
-g(\varphi_{1}(\rho {X}),\overline{Y})\,\overline{u}-\textbf{T}(\overline{u},\rho {X},\overline{Y})\,\varphi_{1}(\overline{\eta}) \nonumber .
\end{eqnarray*}

\item{(8)} If $f_{1}=1$,  $A=u$ and $\varphi_{2}=0$, then we get a Finslerian version of the connection $\overline{D}$ presented in \cite[\S 4.2 (5)]{New}. Therefore, $\overline{D}$ is given by:
 \begin{eqnarray*}
\overline{D} _{X} \overline{Y}&=&\nabla _{X} \overline{Y}+g(\rho {X}-\varphi_{1}(\rho {X}),\overline{Y})\,\overline{u} -u(\rho {X})\overline{Y}-u(\overline{Y})\,\rho {X}-L\,\ell(\overline{Y})\,T(\overline{u},\rho {X})\nonumber\\
&&+\textbf{T}(\overline{u},\rho {X},\overline{Y})\overline{\eta}
+L^{2}\,S(\rho{X},\overline{u})\overline{Y} + L\,\ell(\varphi_{1}(\overline{Y}))T(\overline{u},\rho {X})
\nonumber\\
&&
-u(\overline{\eta})\{S(\varphi_{1}(\overline{\eta}),\rho {X} ) \overline{Y} +T(\varphi_{1}(\overline{Y}),\rho {X}) -\varphi_{1}(T(\rho {X},\overline{Y}))\}
\nonumber\\
&& +u(\overline{Y})\,\varphi_{1}(\rho {X}) +L\,\ell(\varphi_{1}(\overline{\eta}))\,S(\overline{u},\rho {X})\overline{Y}
-\textbf{T}(\overline{u},\rho {X},\overline{Y})\,\varphi_{1}(\overline{\eta}) \nonumber .
\end{eqnarray*}
\item{(9)} If $\varphi_{1}=0$, then we obtain a Finslerian version of the quarter-symmetric recurrent connection $\overline{D}$ given in \cite[\S 4.2 (6)]{New}. Thus, $\overline{D}$ has the form
\begin{eqnarray*}
\overline{D} _{X} \overline{Y}&=&\nabla _{X} \overline{Y}
+f_{1}\{g(\rho {X},\overline{Y})\overline{a}-A(\rho {X})\overline{Y}-A(\overline{Y})\,\rho {X}-L\ell(\overline{Y})T(\overline{a},\rho {X})\nonumber\\
&&+\textbf{T}(\overline{a},\rho {X},\overline{Y})\overline{\eta}
+L^{2}\,S(\rho {X},\overline{a})\overline{Y}\}-u(\rho {X})\varphi_{2}(\overline{Y})+u(\overline{Y})T(\varphi_{2}(\overline{\eta}),\overline{X})
\nonumber \\
&& -\mathbf{T}(\varphi_{2}(\overline{\eta}),\rho {X},\overline{Y})\overline{u}
+u(\overline{\eta})\,S(\varphi_{2}(\overline{\eta}),\rho {X})\overline{Y}.
\end{eqnarray*}
\item{(10)} When  $f_{1}=1$, $A=u$ and $\varphi_{1}=0,\,\,\overline{D}$ reduces to a Finslerian version of the special quarter-symmetric recurrent connection presented  in  \cite[\S 4.2 (7)]{New} . That is,
\begin{eqnarray*}
\overline{D} _{X} \overline{Y}&=&\nabla _{X} \overline{Y}
+g(\rho {X},\overline{Y})\overline{u}-u(\rho {X})\overline{Y}-u(\overline{Y})\rho {X}-L\ell(\overline{Y})T(\overline{u},\rho {X})\nonumber\\
&&+\textbf{T}(\overline{u},\rho {X},\overline{Y})\overline{\eta}
+L^{2}\,S(\rho {X},\overline{u})\overline{Y}\nonumber -u(\rho {X})\varphi_{2}(\overline{Y})+u(\overline{Y})T(\varphi_{2}(\overline{\eta}),\overline{X})
\\
&&
 -\mathbf{T}(\varphi_{2}(\overline{\eta}),\rho {X},\overline{Y})\overline{u}
+u(\overline{\eta})\,S(\varphi_{2}(\overline{\eta}),\rho {X})\overline{Y}.
\end{eqnarray*}
\end{description}
\newpage
\begin{flushleft}
\textbf{Quarter-symmetric non-metric Finsler connection: $f_{1}=0 ,\,f_{2}\neq 0$}
\end{flushleft}
\begin{description}
\item{(11)} If  $\varphi_{2}=0$,  we get a Finslerian version of the connection given in \cite[\S 4.2 (8)]{New}. Then, the connection $\overline{D}$ has the form:
\begin{eqnarray*}
\overline{D} _{X} \overline{Y}&=&\nabla _{X} \overline{Y} -f_{2}\{g(\rho {X},\overline{Y})\,\overline{b}-L\,\ell(\overline{Y})\,T(\overline{b},\rho {X})
+\textbf{T}(\overline{b},\rho {X},\overline{Y})\,\overline{\eta} \\
&&+L^{2}\,S(\rho{X},\overline{a})\overline{Y} \}-g(\varphi_{1}(\rho {X}),\overline{Y})\overline{u}
+ L\,\ell(\varphi_{1}(\overline{Y}))T(\overline{u},\rho {X})\\
&&-u(\overline{\eta})\{S(\varphi_{1}(\overline{\eta}),\rho {X} ) \overline{Y} +T(\varphi_{1}(\overline{Y}),\rho {X})
  -\varphi_{1}(T(\rho {X},\overline{Y}))\}\\
&& +u(\overline{Y}\,\varphi_{1}(\rho {X}) +L\,\ell(\varphi_{1}(\overline{\eta}))\,S(\overline{u},\rho {X})\overline{Y}
-\textbf{T}(\overline{u},\rho {X},\overline{Y})\,\varphi_{1}(\overline{\eta}) .
\end{eqnarray*}
\item{(12)} When $B=u$  and $\varphi_{2}=0,\,\,\overline{D}$  becomes a Finslerian version of the connection appeared in  \cite[\S 4.2 (9)]{New}. That is,
\begin{eqnarray*}
\overline{D} _{X} \overline{Y}&=&\nabla _{X} \overline{Y} -f_{2}\{g(\rho {X},\overline{Y})\,\overline{u}-L\,\ell(\overline{Y})\,T(\overline{u},\rho {X})
+\textbf{T}(\overline{u},\rho {X},\overline{Y})\,\overline{\eta} \\
&&+L^{2}\,S(\rho{X},\overline{a})\overline{Y} \}-g(\varphi_{1}(\rho {X}),\overline{Y})\overline{u}
+ L\,\ell(\varphi_{1}(\overline{Y}))T(\overline{u},\rho {X})\\
&&-u(\overline{\eta})\{S(\varphi_{1}(\overline{\eta}),\rho {X} ) \overline{Y} +T(\varphi_{1}(\overline{Y}),\rho {X})
  -\varphi_{1}(T(\rho {X},\overline{Y}))\}\\
&& +u(\overline{Y})\,\varphi_{1}(\rho {X})+L\,\ell(\varphi_{1}(\overline{\eta}))\,S(\overline{u},\rho {X})\overline{Y}
-\textbf{T}(\overline{u},\rho {X},\overline{Y})\,\varphi_{1}(\overline{\eta}).
\end{eqnarray*}


\item{(13)} If  $\varphi_{1}=0$,  we obtain a Finslerian version of the connection given in \cite[\S 4.2 (10)]{New}.
 The connection $\overline{D}$ can be written as
\begin{eqnarray*}
\overline{D} _{X} \overline{Y}&=&\nabla _{X} \overline{Y} -f_{2}\{g(\rho {X},\overline{Y})\,\overline{b}-L\,\ell(\overline{Y})\,T(\overline{b},\rho {X})+\textbf{T}(\overline{b},\rho {X},\overline{Y})\,\overline{\eta} \\&&+L^{2}\,S(\rho{X},\overline{a})\overline{Y} \}-\mathbf{T}(\varphi_{2}(\overline{\eta}),\rho {X},\overline{Y})\,\overline{u}-u(\rho {X})\,\varphi_{2}(\overline{Y})\\&&+u(\overline{\eta})\,S(\varphi_{2} (\overline{\eta}),\rho {X} ) \overline{Y}  +u(\overline{Y})\,T(\varphi_{2}(\overline{\eta}),\overline{X}) .
\end{eqnarray*}

\item{(14)} If $B=u$ and $\varphi_{1}=0$, we get a Finslerian version of the connection presented in \cite[\S 4.2 (11)]{New}. The connection $\overline{D}$ is given by:
\begin{eqnarray*}
\overline{D} _{X} \overline{Y}&=&\nabla _{X} \overline{Y} -f_{2}\{g(\rho {X},\overline{Y})\,\overline{u}-L\,\ell(\overline{Y})\,T(\overline{u},\rho {X})+\textbf{T}(\overline{u},\rho {X},\overline{Y})\,\overline{\eta} \\&&+L^{2}\,S(\rho{X},\overline{a})\overline{Y} \}-\mathbf{T}(\varphi_{2}(\overline{\eta}),\rho {X},\overline{Y})\,\overline{u}-u(\rho {X})\,\varphi_{2}(\overline{Y})\\&&+u(\overline{\eta})\,S(\varphi_{2} (\overline{\eta}),\rho {X} ) \overline{Y}  +u(\overline{Y})\,T(\varphi_{2}(\overline{\eta}),\overline{X}) .
\end{eqnarray*}
\end{description}
\begin{flushleft}
\textbf{Semi-symmetric metric Finsler connection: $f_{1}=f_{2}=0,\,\varphi=id_{\pi^{-1}(TM)}$}
\end{flushleft}
\begin{description}
\item{(15)} We obtain the Finslerian version of the connection defined in   \cite{Yano.2}.
 That is,\vspace{-0.1 cm}
\begin{eqnarray*}
\overline{D} _{X} \overline{Y}&=&\nabla _{X} \overline{Y}-g(\rho {X},\overline{Y})\overline{u}
+u(\overline{Y})\rho {X}+L\ell(\overline{Y})T(\overline{u},\rho {X})\nonumber\\
&&+L^{2}T(T(\rho {X},\overline{Y}),\overline{u})-L^{2}T(T(\overline{u},\overline{Y}),\rho {X})
-T(\overline{u},\rho {X},\overline{Y})\overline{\eta}.
\end{eqnarray*}

\item{(16)} If $u=\ell$, we get an intrinsic version of the connection given in \cite[\S 5.3]{amr20}.
 \vspace{-0.3 cm}
\begin{equation*}
\overline{D} _{X} \overline{Y}=\nabla _{X} \overline{Y}-L^{-1}g(\rho {X},\overline{Y})\overline{\eta}
+\ell(\overline{Y})\rho {X}.
\end{equation*}
\end{description}

\begin{flushleft}
\textbf{Semi-symmetric non-metric Finsler connection: $f_{1}\neq0,\,f_{2}=0,\varphi=id_{\pi^{-1}(TM)}$}
\end{flushleft}
\begin{description}
\item{(17)}   $\overline{D}$ reduces to a Finslerian version of the semi-symmetric recurrent connection given in \cite{amr20, New}.  Thus, $\overline{D}$ is  written in the form:
\begin{eqnarray*}
\overline{D} _{X} \overline{Y}&=&\nabla _{X} \overline{Y}+f_{1}\{g(\rho {X},\overline{Y})\overline{a}-A(\rho {X})\overline{Y}-A(\overline{Y})\,\rho {X}-L\ell(\overline{Y})T(\overline{a},\rho {X})\nonumber\\
&&+\textbf{T}(\overline{a},\rho {X},\overline{Y})\overline{\eta}
+L^{2}\{{T}(T(\overline{a},\overline{Y}),\rho {X})-{T}(T(\rho {X},\overline{Y}),\overline{a})\}\}\nonumber \\
&&-g(\rho {X},\overline{Y})\overline{u}+u(\overline{Y})\rho {X}
+L\ell(\overline{Y})T(\overline{u},\rho {X}) \nonumber\\
&&-L^{2}T(T(\overline{u},\overline{Y}),\rho {X})-\textbf{T}(\overline{u},\rho {X},\overline{Y})\overline{\eta}
 +L^{2}T(T(\rho {X},\overline{Y}),\overline{u}).
\end{eqnarray*}

\item{(18)} If $f_{1}=\frac{1}{2}$,  we obtain a Finslerian version of  the semi-symmetric recurrent connection $\overline{D}$
studied in \cite{Andonie, Liang, amr20}.
 In this case,  $\overline{D}$ is given by:
\begin{eqnarray*}
\overline{D} _{X} \overline{Y}&=&\nabla _{X} \overline{Y}+\frac{1}{2}\{g(\rho {X},\overline{Y})\overline{a}-A(\rho {X})\overline{Y}-A(\overline{Y})\,\rho {X}-L\ell(\overline{Y})T(\overline{a},\rho {X})\nonumber\\
&&+\textbf{T}(\overline{a},\rho {X},\overline{Y})\overline{\eta}
+L^{2}\{{T}(T(\overline{a},\overline{Y}),\rho {X})-{T}(T(\rho {X},\overline{Y}),\overline{a})\}\}\nonumber \\
&&-g(\rho {X},\overline{Y})\overline{u}+u(\overline{Y})\rho {X}+ L\ell(\overline{Y})T(\overline{u},\rho {X})\nonumber\\
&&+L^{2}T(T(\rho {X},\overline{Y}),\overline{u}) -L^{2}T(T(\overline{u},\overline{Y}),\rho {X})-T(\overline{u},\rho {X},\overline{Y})\overline{\eta}.
\end{eqnarray*}

\item{(19)} If $f_{1}=\frac{1}{2}$ and $A=u=\ell$,   we obtain a special semi-symmetric $h$-recurrent Finsler connection \cite{amr20}
 given by:
\[\overline{D} _{X} \overline{Y}=\nabla _{X} \overline{Y}-\frac{1}{2}\{L^{-1}g(\rho {X},\overline{Y})\overline{\eta}+\ell(\rho {X})\overline{Y}-\ell(\overline{Y})\rho {X}\}.\]
\end{description}

\begin{flushleft}
\textbf{Semi-symmetric non-metric Finsler connection:$\,\,f_{1}=0 ,\,f_{2}\neq0, \,\varphi=id_{\pi^{-1}(TM)}$}
\end{flushleft}
\begin{description}
\item{(20)}  We obtain a Finslerian version of the semi-symmetric non-metric connection given in \cite[\S 4.4 (14)]{New},  that is,
\begin{eqnarray*}
\overline{D} _{X} \overline{Y}&=&\nabla _{X} \overline{Y}-f_{2}\{g(\rho {X},\overline{Y})\overline{b}-L\ell(\overline{Y})T(\overline{b},\rho {X})
+\textbf{T}(\overline{b},\rho {X},\overline{Y})\overline{\eta}\nonumber \\
&&+L^{2}\{{T}(T(\overline{b},\overline{Y}),\rho {X})-{T}(T(\rho {X},\overline{Y}),\overline{b})\}\}\nonumber \\
&&-g(\rho {X},\overline{Y})\overline{u}+u(\overline{Y})\rho {X}+ L\ell(\overline{Y})T(\overline{u},\rho {X})\nonumber\\
&&-L^{2}T(T(\overline{u},\overline{Y}),\rho {X})-\textbf{T}(\overline{u},\rho {X},\overline{Y})\overline{\eta}
+L^{2}T(T(\rho {X},\overline{Y}),\overline{u}).
\end{eqnarray*}

\item{(21)} If $f_{2}=-1$,   we obtain a Finslerian version of the connection discussed in \cite{Sengupta},  that is,
\begin{eqnarray*}
\overline{D} _{X} \overline{Y}&=&\nabla _{X} \overline{Y}+ g(\rho {X},\overline{Y})\overline{b}-L\ell(\overline{Y})T(\overline{b},\rho {X})
+\textbf{T}(\overline{b},\rho {X},\overline{Y})\overline{\eta}\nonumber \\
&&+L^{2}\{{T}(T(\overline{b},\overline{Y}),\rho {X})-{T}(T(\rho {X},\overline{Y}),\overline{b})\}\nonumber \\
&&-g(\rho {X},\overline{Y})\overline{u}+u(\overline{Y})\rho {X}+ L\ell(\overline{Y})T(\overline{u},\rho {X})\nonumber\\
&&-L^{2}T(T(\overline{u},\overline{Y}),\rho {X})-\textbf{T}(\overline{u},\rho {X},\overline{Y})\overline{\eta}
+L^{2}T(T(\rho {X},\overline{Y}),\overline{u}).
\end{eqnarray*}

\item{(22)} If $f_{2}=-1$ and $B=u$, then we get the Finslerian version of the connection studied in  \cite{Agashe}, that is,
\begin{eqnarray*}
\overline{D} _{X} \overline{Y}&=&\nabla _{X} \overline{Y}+u(\overline{Y})\rho {X}.
\end{eqnarray*}
\end{description}

\begin{flushleft}
\textbf{Symmetric non-metric Finsler connection: $u=0$}
\end{flushleft}
\begin{description}
\item{(23)} We obtain a Finslerian version of  the connection appeared in \cite[\S 4.5 (15)]{New}. It is given by: \vspace{-0.2cm}
\begin{eqnarray*}
\overline{D} _{X} \overline{Y}&=&\nabla _{X} \overline{Y}+f_{1} \,\{g(\rho {X},\overline{Y})\,\overline{a}-A(\rho {X})\overline{Y}-A(\overline{Y})\,\rho {X}-L\,\ell(\overline{Y})\,T(\overline{a},\rho {X})\nonumber\\
&&+\textbf{T}(\overline{a},\rho {X},\overline{Y})\overline{\eta}
+L^{2}\,S(\rho{X},\overline{a})\overline{Y} \} -f_{2}\{g(\rho {X},\overline{Y})\,\overline{b}-L\,\ell(\overline{Y})\,T(\overline{b},\rho {X})
\nonumber\\
&&+\textbf{T}(\overline{b},\rho {X},\overline{Y})\,\overline{\eta}+L^{2}\,S(\rho{X},\overline{b})\overline{Y} \}  \nonumber.
\end{eqnarray*}

\item{(24)} If $f_{1}=\frac{1}{2}$  and $f_{2}=0$, then $\overline{D}$ reduced to a Finslerian version of the symmetric recurrent connection (Weyl
 connection)
investigated in  \cite{Folland, hv, h-rec.}.  That is, \vspace{-0.2cm}
\begin{eqnarray*}
\overline{D} _{X} \overline{Y}&=&\nabla _{X} \overline{Y}+\frac{1}{2}\{g(\rho {X},\overline{Y})\,\overline{a}-A(\rho {X})\,\overline{Y}-A(\overline{Y})\,\rho {X}-L\,\ell(\overline{Y})\,T(\overline{a},\rho {X})\nonumber\\
&&\qquad \qquad \quad+\textbf{T}(\overline{a},\rho {X},\overline{Y})\,\overline{\eta}
+L^{2}\, S(\rho {X},\overline{a})\overline{Y}\}.
\end{eqnarray*}

\item{(25)} If $f_{1}=f_{2}=-1$  and $A=B$,  then $\overline{D}$ is a Finslerian version of  the connection  considered in \cite{Yano.3}, that is,\vspace{-0.1cm}
\begin{equation*}
\overline{D}_{X} \overline{Y}=\nabla _{X} \overline{Y}+A(\rho {X})\overline{Y}+A(\overline{Y})\,\rho {X}.
\end{equation*}

\item{(26)} If $f_{1}=\frac{1}{2}$, $A=\ell$  and $f_{2}=0$,  then we obtain a special symmetric h-recurrent Finsler connection studied in  \cite{hv, h-rec.}, that is,\vspace{-0.2cm}
\begin{equation*}
\overline{D} _{X} \overline{Y}=\nabla _{X} \overline{Y}+\frac{1}{2}\{L^{-1}g(\rho {X},\overline{Y})\overline{\eta}-\ell(\rho {X})\overline{Y}-\ell(\overline{Y})\rho {X}\}.
\end{equation*}
\end{description}
We end this work by the following remark: applying the ${P}^{1}$-process and ${C}$-process to each of the above mentioned special cases,  one can get more new Finsler connections.
\begin{flushleft}
\textbf{Acknowledgment} We would like to express our deep thanks to Professor Nabil L. Youssef for the careful reading of the manuscript and his valuable suggestions which led to this version.
\end{flushleft}
\providecommand{\bysame}{\leavevmode\hbox
to3em{\hrulefill}\thinspace}
\providecommand{\MR}{\relax\ifhmode\unskip\space\fi MR }
\providecommand{\MRhref}[2]{%
  \href{http://www.ams.org/mathscinet-getitem?mr=#1}{#2}
} \providecommand{\href}[2]{#2}

\end{document}